\documentclass[12pt]{article}

\usepackage{amstext    }
\usepackage{amsthm    }
\usepackage{a4}
\usepackage[mathscr]{eucal} 
\usepackage{mathrsfs}

\usepackage{amsmath}
\usepackage{amssymb}
\usepackage{amscd}

\numberwithin{equation}{section}

\newtheorem{theorem}{Theorem}[section]
\newtheorem{definition}[theorem]{Definition}
\newtheorem{proposition}[theorem]{Proposition}
\newtheorem{corollary}[theorem]{Corollary}
\newtheorem{lemma}[theorem]{Lemma}
\newtheorem{remark}[theorem]{Remark}

\newcommand{\cali}[1]{\mathscr{#1}}

\newcommand{\supp}{{\rm supp}}
\newcommand{\const}{{\rm const}}

\newcommand{\loc}{{loc}}
\newcommand{\ddc}{dd^c}

\newcommand{\DSH}{{\rm DSH}}

\newcommand{\BB}{{\rm B}}

\newcommand{\Cc}{\cali{C}}

\newcommand{\Lc}{\cali{L}}

\newcommand{\C}{\mathbb{C}}

\newcommand{\R}{\mathbb{R}}

\title{Pull-back of currents by holomorphic maps}
\author{Tien-Cuong Dinh and Nessim Sibony}

\begin{document}

\maketitle

\begin{abstract}
We define the pull-back operator, associated to a meromorphic transform,
on several types of currents. We also give a simple proof to a 
version of a classical theorem on the extension of currents.
\end{abstract}


%%%%%%%%%%%%%%%%%%%%%%%%%%%%%%%%%%%%%%%%%%%%%%%%%%%%%

\section{Introduction} \label{section_introduction}

Let $X$ and $X'$ be two (connected) complex
manifolds of dimensions $k$ and $k'$ respectively. 
A holomorphic map $f:X\rightarrow X'$ induces a pull-back
operator $f^*$ acting on smooth forms on $X'$ with values in
the space of smooth forms on $X$. 
Our main purpose in this paper is to extend the previous operator to
some classes of currents. This is a  fundamental question in complex
analysis because in many problems, one has to deal with singular
objects like subvarieties or more generally 
positive closed currents.

Our first motivation comes from the theory of complex dynamical
systems in higher dimension. Given a holomorphic self-map
or more generally a multivalued meromorphic self-map $f:X\rightarrow
X$, a main problem in dynamics is
to construct interesting measures invariant by $f$. A
general strategy is to construct invariant positive closed currents
and then obtain invariant measures as intersection of such
currents, see \cite{Sibony1} for historical comments. 
So, a necessary step here is to define the pull-back operator
on positive closed currents. Theorem \ref{th_main} below allows to
extend the contruction of Green currents in \cite{DinhSibony4} to arbitrary
holomorphic correspondences with finite fibers between compact K\"ahler manifolds.

In order to avoid some trivial counter-examples, assume throughout the paper 
that $f$ is dominant, i.e. its image contains a non-empty open subset
of $X'$. Otherwise, one might replace $X'$ by $f(X)$ which is a
variety 
immersed in $X'$, possibly with singularities. In
\cite{Meo}, M\'eo proved that the operator $f^*$ can be continuously  extended
to positive closed currents of bidegree $(1,1)$. He also
constructed an example where one cannot extend $f^*$ to
positive closed currents of higher bidegree. 
Our main theorem is the following result.

\begin{theorem} \label{th_main}
Let $f:X\rightarrow X'$ be a holomorphic map. Assume that each
fiber of $f$ is either empty or is an analytic set of dimension $\dim
X-\dim X'$. Then the pull-back operator $f^*$ can be extended to
positive closed (resp. $\ddc$-closed) $(p,p)$-currents on $X'$. Moreover, if $T$ is such a current
then $f^*(T)$ is a positive closed (resp. $\ddc$-closed) $(p,p)$-current on $X$ which
depends continuously on $T$ for the weak topology on currents. If $T$
has no mass on a Borel set $A\subset X'$, 
then $f^*(T)$ has no mass on $f^{-1}(A)$. 
\end{theorem} 

Observe that $f$ is open hence we can consider that $f$ is surjective
by restricting $X'$.
Note that when $f$ is a finite map between open subsets of $\C^k$, M\'eo gave in
\cite{Meo} a definition of $f^*$ for positive closed $(p,p)$-currents. 
He used potentials of currents and didn't consider the continuity
of $f^*$ and its independence of coordinates in $\C^k$ which are
crucial here in order to extend $f^*$ to the case of manifolds. 

We deduce from the previous result the following corollary  where 
$\{\cdot\}$ denotes the class of a positive $\ddc$-closed $(p,p)$-current in
the cohomology group $H^{p,p}(\cdot, \C)$:

\begin{corollary} \label{cor_main}
Let $f:X\rightarrow X'$ be a holomorphic map between
  compact K\"ahler manifolds. Assume that each
fiber of $f$ is an analytic set of dimension $\dim
X-\dim X'$.  If $T$ is a positive closed or
  $\ddc$-closed 
  $(p,p)$-current on $X'$, then $f^*\{T\}=\{f^*(T)\}$.
\end{corollary}
 
The linear operator $f^*:H^{p,p}(X',\C)\rightarrow H^{p,p}(X,\C)$ is
induced by the action of $f^*$ on smooth forms. Corollary \ref{cor_main} follows from the
continuity in Theorem \ref{th_main} and from our 
result in \cite{DinhSibony3} which says that one can write $T=T^+-T^-$
where $T^\pm$ can be approximated
by smooth positive $\ddc$-closed forms.

Theorem \ref{th_main} is also valid for
another class of currents, useful in dynamics, that we call {\it dsh
currents}, see Theorem \ref{th_pullback_dsh}. 
It is still valid for some meromorphic maps or more generally for
some meromorphic transforms, see Theorems \ref{th_pullback_pure} and
\ref{th_pullback_pure_dsh}. 
The case of compact K\"ahler manifolds will be discussed in Section \ref{section_karler}.

The tools in order to prove the main results give also a simple proof of the 
following theorem on the extension of currents.

\begin{theorem} \label{th_extension}
Let $F$ be a
closed subset of a complex manifold $X$. Let $T$ be a positive
$(p,p)$-current on $X\setminus F$.
Assume that $F$ is  locally complete pluripolar and 
$T$ has locally finite mass near $F$. Assume also that there exists a positive 
$(p+1,p+1)$-current $S$ with locally finite mass near $F$ such that $\ddc T\leq S$ on $X\setminus
F$.
Then $\ddc T$ has locally finite mass near $F$. Moreover, if 
$\widetilde T$ and $\widetilde{\ddc T}$ denote
the extensions by zero of $T$ and $\ddc T$ on $X$, then  
$\widetilde{\ddc T} -\ddc \widetilde T$ is positive. If $T$ is closed then $\widetilde T$ is closed.
\end{theorem}
This result extends a classical Skoda's theorem when $T$ is closed,
see \cite{Skoda1}, \cite{ElMir} and \cite{Sibony2}.   
For $S= 0$ it is proved by Dabbek-Elkhadhara-El Mir in
\cite{DabbekElkhadharaElMir}, see also Remark \ref{rk_extension}, and is due to
Alessandrini-Bassanelli \cite{AlessandriniBassanelli1, Bassanelli} when  $F$ is an
analytic set and $\ddc T$ has bounded mass. 
Under the extra assumption that $dT$ is of order zero, the result 
was proved by the second author in \cite{Sibony2}.
In this case, we also have the formula $d\widetilde T= \widetilde{dT}$.

%%%%%%%%%%%%%%%%%%%%%%%%%%%%%%%%%%%%%%%%%%%%%%

\section{Extension of currents}

In this section we give a simple proof of Theorem
\ref{th_extension}. 
We start with the following lemmas which are versions of the
Chern-Levine-Nirenberg inequality \cite{ChernLevineNirenberg,
  Demailly2}.
In what follows, $\omega$ denotes a hermitian $(1,1)$-form on a
manifold $X$ of dimension $k$. 
If $T$ is a current of order zero, the mass of $T$ on a Borel set
$K\subset X$ is denoted by $\|T\|_K$. The mass of $T$ on
$X$ is denoted by $\|T\|$.
If $T$ is a positive or a negative 
$(p,p)$-current, $\|T\|_K$ is equivalent to
$|\int_K T\wedge \omega^{k-p}|$. We often identify these two quantities.

\begin{lemma} \label{lemma_cln}
Let $U$ be an open subset of $X$.
Let $K$ and $L$ be compact sets in $U$ with $L\Subset K$. 
Assume that $T$ is positive and $\ddc T$ has order zero. 
Then there exists  a constant $c_{K,L}>0$ such that
for every smooth bounded psh function $u$ on $U$ we have the following estimate
$$\int_L du \wedge d^c u \wedge T\wedge \omega^{k-p-1}
\leq c_{K,L} \|u\|_{\Lc^\infty(K)}^2\big( \|T\|_K + \|\ddc T\|_K \big).$$
\end{lemma}
\proof
We can assume that $K$ is the unit ball in $\C^k$
and $L$ is the ball of center 0 and of radius $1-3\delta$,
$0<\delta<1/4$, and $\omega$ is the canonical K\"ahler form on $\C^k$.
Replacing $T$ by $T\wedge \omega^{k-p-1}$ and $u$ by ${1\over
  4}\|u\|_{\Lc^\infty(K)}^{-1}u+{1\over 4}$,
we can assume that $p=k-1$ and $0\leq u\leq 1/2$. 
 Let $\chi$ be a real-valued smooth function, 
$0\leq \chi\leq 1$, supported on $K$ and such that $\chi=1$ on $\{\|z\|<1-\delta\}$. 
Define 
$$v(z):={\max}_\epsilon\big(u(z), \delta^{-1}(\|z\|^2-(1-2\delta)^2)\big)$$
where ${\max}_\epsilon(x,y)$ is a smooth function on $\R^2$, convex increasing on
each variable, $\max_\epsilon(x,y)\geq \max(x,y)$, and
$\max_\epsilon(x,y)$ is equal to $\max(x,y)$ outside a small neighbourhood of $\{x=y\}$.
Then $v$ is a positive 
smooth psh function equal to $\delta^{-1}(\|z\|^2-(1-2\delta)^2)$ on $\{1-\delta\leq \|z\|\leq 1\}$.
In particular, we have $v(z)=\delta^{-1}(\|z\|^2-(1-2\delta)^2)$ on
the support $\supp(d\chi)$ of $d\chi$ and $v=u$ on $L$. Then, since $v^2$ is
psh, 
\begin{eqnarray*}
\int_L du\wedge d^c u \wedge T
 & = & \int_L \Big( {1\over 2} \ddc v^2  -v\ddc v\Big) \wedge T \nonumber\\
& \leq & {1 \over 2}\int \ddc v^2 \wedge\chi T 
 =   {1 \over 2} \int v^2 \ddc (\chi T) \nonumber\\
& = & {1 \over 2} \int  v^2 (\chi\ddc T + T\wedge \ddc \chi) 
 + {1\over 2} \int v^2 (-d^c\chi\wedge dT - d^c T\wedge d\chi).
\end{eqnarray*}
Define $\widetilde v(z):=\delta^{-1}(\|z\|^2-(1-2\delta)^2)$ on $\C^k$.
Recall that 
$v=\widetilde v$ on $\supp(d\chi)$. 
Using an analogous computation as above, we obtain that the last integral is equal to
$$\int \ddc \widetilde v^2 \wedge \chi T - 
\int \chi \widetilde v^2 \ddc T - \int \widetilde v^2 \ddc \chi \wedge T.$$
Hence
$$\int_L du\wedge d^c u \wedge T\leq {1 \over 2} \int  (v^2-\widetilde v^2) (\chi\ddc T + T\wedge \ddc \chi) 
 + {1\over 2}\int \ddc \widetilde v^2 \wedge \chi T.$$
The lemma follows.
\endproof

\begin{lemma} \label{lemma_cln_bis}
Under the assumptions of Lemma \ref{lemma_cln}, we have
$$\|\ddc u\wedge T\|_L\leq c_{K,L}\|u\|_{\Lc^\infty(K)} \big( \|T\|_K + \|\ddc T\|_K \big)$$
where $c_{K,L}>0$ is a constant independent of $u$ and $T$. 
\end{lemma}
\proof
Write 
\begin{eqnarray}
\ddc u\wedge T := u\ddc T -\ddc (uT) +d(d^c u\wedge T) - d^c(du\wedge T). \label{eq_ddc}
\end{eqnarray}
Hence, Cauchy-Schwarz inequality and Lemma \ref{lemma_cln}, applied to $L':=\supp(\chi)$,  imply
\begin{eqnarray*}
\|\ddc u\wedge T\|_L &\leq & \int \chi \ddc u\wedge T \\
& = &  \int \chi u \ddc T -\int \ddc \chi \wedge uT
- \int d \chi \wedge d^c u\wedge T + \int d^c \chi \wedge  du \wedge T \\
& \lesssim & \|u\|_{\Lc^\infty(K)} \big( \|T\|_K + \|\ddc T\|_K \big).
\end{eqnarray*}
\endproof

\noindent
{\it End of the proof of Theorem \ref{th_extension}---} 
Recall that a {\it complete locally pluripolar set} is locally the pole set
$\{\varphi=-\infty\}$ of a psh function $\varphi$.
Since the problem is local, we can assume that $X$ is a ball in $\C^k$. 
Then, there is a sequence of smooth psh functions $u_n$,
$0\leq u_n\leq 1$, vanishing near $F$ and increasing to 1 on
$X\setminus F$, see e.g. \cite{Sibony1}.
We have $u_n T\rightarrow \widetilde T$  and as in (\ref{eq_ddc})
\begin{eqnarray}
u_n\ddc T-\ddc(u_nT)  = \ddc u_n\wedge T -d (d^c u_n\wedge T) + d^c ( du_n\wedge T).
\label{eq_difference_ddc}
\end{eqnarray}
We first show that $\partial u_n\wedge T\rightarrow 0$.
This implies that $\overline\partial u_n\wedge T\rightarrow 0$ by
conjugaison. We can assume that $p=k-1$. Indeed, one only has to test
forms of type $\varphi\Omega$ with $\varphi$ a smooth function and
$\Omega$ a smooth positive closed form. Then we reduce the problem to
the case $p=k-1$ by replacing $T$ by $T\wedge \Omega$.

Let $\alpha$ be a test smooth $(0,1)$-form with support in a compact
subset $K$ of $X$. 
We have to show that $\int\partial u_n\wedge\alpha\wedge T$ converge
to zero.
Fix an $\epsilon>0$ and
a neighbourhood $U_\epsilon$ of $F$ such that 
$$\int_{U_\epsilon} i\alpha\wedge \overline\alpha\wedge \widetilde T\leq \epsilon^2.$$
Since $u_n$ converge locally uniformly to 1 on $X\setminus F$, Lemma \ref{lemma_cln}, applied to $(u_n-1)$, 
implies that locally on $X\setminus F$ the mass of 
$\partial u_n \wedge T$ tends to zero.
Hence,
$$\lim_{n\rightarrow\infty} \int_{X\setminus U_\epsilon}  \partial u_n \wedge\alpha\wedge T =0.$$
On the other hand, Cauchy-Schwarz inequality gives
 $$\Big|\int_{U_\epsilon} i\partial u_n \wedge \alpha \wedge  T \Big|\leq
\Big(\int_K i\partial u_n \wedge \overline\partial u_n \wedge T \Big)^{1/2}
\Big(\int_{U_\epsilon} i\alpha \wedge \overline \alpha \wedge  T \Big)^{1/2}.$$
The second factor 
is bounded by $\epsilon$. So it is enough to show that the integral on $K$ is uniformly bounded with respect to $n$.
We cannot apply directly Lemma \ref{lemma_cln} since we don't know that $\ddc
T$ has finite mass.
Let $\chi$ be a cutoff function, $0\leq \chi\leq 1$ and $\chi=1$ on $K$. 
We have, using a computation as in (\ref{eq_ddc}),
\begin{eqnarray*}
\lefteqn{I_n:=\int\chi^2  i\partial u_n \wedge \overline\partial u_n \wedge T}\\
 & \leq & {1\over 2} \int\chi^2 
i\partial\overline\partial u_n^2\wedge T = {1\over 2} \int u_n^2 i\partial\overline\partial (\chi^2 T)\\
& \leq &  {1\over 2} \int u_n^2 \big[\chi^2 i\partial\overline\partial T + i\partial \overline \partial \chi^2
\wedge T \big] + 2\Big|\int u_n\overline \partial u_n \wedge \partial\chi^2 \wedge T\Big|
\end{eqnarray*}
The first term in the last line is bounded uniformly since $0\leq u_n\leq 1$ and $\ddc T\leq S$. 
The identity $\partial \chi^2=
2\chi\partial\chi$ and the Cauchy-Schwarz inequality imply that the last integral is bounded by
$$2\Big(\int\chi^2 i\partial u_n\wedge\overline\partial u_n\wedge T\Big)^{1/2} \Big( \int i\partial\chi \wedge 
\overline\partial\chi \wedge u_n^2 T\Big)^{1/2}.$$
Since the last integral is uniformly bounded we have $I_n\leq \const (1+I_n^{1/2})$ and hence $I_n$ is bounded.

We have proved that $\partial u_n\wedge T\rightarrow 0$ and $\overline\partial u_n\wedge T\rightarrow 0$.
Identity (\ref{eq_difference_ddc}) implies that 
$$u_n(\ddc T-S) - \ddc u_n\wedge T=-u_nS+\ddc(u_nT)-d(d^cu_n\wedge T)+d^c(du_n\wedge T).$$
The right hand side converges to $-\widetilde S+\ddc \widetilde T$ where $\widetilde S$ is the trivial extension
by zero of $S$ on $X$. Since both terms on the 
left hand side are negative currents their limit values are
negative. We then deduce that
$-\widetilde S+\ddc \widetilde T$ is negative and hence
$\ddc T$ has finite mass near $F$.
Finally the left hand side of (\ref{eq_difference_ddc}) converges to
$\widetilde{\ddc T}-\ddc \widetilde T$ and hence the 
right hand side converge to $\lim \ddc u_n\wedge T$ which is positive. 
Hence $\widetilde{\ddc T} - \ddc \widetilde T$ is positive.
When $T$ is closed, we have $d\widetilde T= \lim d u_n\wedge T=0$, hence $\widetilde T$ is closed. 
\hfill $\square$

\begin{remark} \label{rk_extension}
\rm
If $S$ is closed then $\widetilde S$ is
closed. It follows from Theorem \ref{th_extension} that $-\ddc \widetilde
T+\widetilde S$ is positive and 
closed. Hence $\widetilde T$ is dsh,
see Section \ref{section_pullback} for the definition.  
If $F$ is an analytic subset of codimension $\geq p+2$ and if $T$ is
$\ddc$-closed then $\widetilde T$ is $\ddc$-closed. Indeed, in this
case, $-\ddc \widetilde T$ is a positive closed $(p+1,p+1)$-current; it
should vanish on sets of codimension $\geq p+2$. 
\end{remark}

%%%%%%%%%%%%%%%%%%%%%%%%%%%%%%%%%%%%%%%%%%%%%%%

\section{Pull-back operator} \label{section_pullback}

We give here the proof of Theorem \ref{th_main}.
We first consider a general local situation.
Let $\BB_k$ denote the unit ball in $\C^k$ and 
$\BB_k(r)$ denote the ball of center 0 and of radius $r$  in $\C^k$.
Let $\pi:\C^k\times \C^{k'}\rightarrow \C^{k'}$ 
be the canonical projection. 
Consider a subvariety $V$ of pure dimension $k'+l$ in $\BB_k\times \BB_{k'}$.
We will denote by $[V]$ the current of integration on $V$.
Assume that the fibers of $\pi_{|V}$ are either empty or of pure dimension $l$.
Let $w=(z,z')$ denote the coordinates in $\C^k\times \C^{k'}$. Let
$\Cc$ denote the critical set of $\pi_{|V}$ which contains the singularities of $V$. 
The set $\Cc$ is defined by the property that 
$\pi_{|V}$ is locally a submersion at every point $w\in V\setminus \Cc$.

\begin{lemma} \label{le_pullback_mass}
Let $L$ be a compact subset of $\BB_k\times \BB_{k'}$. 
Then, there exists $c_L>0$ such that if $T$ is a positive smooth $(p,p)$-form on $\BB_{k'}$ then 
$$\int_{V\cap L}\pi^*(T) \wedge (\ddc\|w\|^2)^{k'+l-p}\leq 
c_L \big(\|T\|+\|\ddc T\|\big).$$
\end{lemma}
\begin{proof}
Since the problem is local we can assume that $0\in V$ and prove the lemma
for a small ball $L$ of center 0.
Let $P$ be a complex plane of codimension $l$ in $\C^k\times\C^{k'}$
such that $\pi$ restricted to $V\cap P$ is discrete.   
Shrinking $\BB_k$ and $\BB_{k'}$ allows to assume that, for every small
pertubation $P_\epsilon$ of $P$, $\pi$ restricted to
$V\cap P_\epsilon$ defines a finite ramified covering of degree $d_t$ over $\BB_{k'}$.
Slicing by $P_\epsilon$ reduces the problem
to the case where $l=0$ and  $\pi_{|V}$ defines a ramified covering of
the same degree $d_t$. 

Observe that 
$u:=(\pi_{|V})_*(\|w\|^2)$ is a continuous psh 
function bounded by $2d_t$. 
Let $K$ be a compact subset of $\BB_{k'}$ such that $L':=\pi(L)\Subset K$.
Then
$$\int_{V\cap L}\pi^*(T) \wedge (\ddc\|w\|^2)^{k'-p} 
\leq \int_{L'} T\wedge  \pi_*\big[(\ddc \|w\|^2)^{k'-p}\big]
\leq \int_{L'} T\wedge  (\ddc u)^{k'-p}.$$ 
Using Lemma \ref{lemma_cln_bis} and the fact that $T$ is smooth, a
simple induction on $r$ gives
$$\int_{L'} T\wedge  (\ddc u)^r\wedge (\ddc\|w\|)^{k'-p-r} \leq
c(\|T\|+\|\ddc T\|),$$
where the constant $c>0$ depends only on $L$ and $d_t$ (this is important
for the slicing). 
When $r=k'-p$, we obtain the lemma.
\end{proof}

We define now the space of dsh currents.

\begin{definition} \label{def_dsh} \rm
A $(p,p)$-current $T$ on a complex manifold $X$ of dimension $k$, 
$0\leq p\leq k-1$, is {\it dsh}
if there exist negative $(p,p)$-currents $T_i$ and positive closed $(p+1,p+1)$-currents $\Omega^\pm_i$
such that 
\begin{eqnarray}
T=T_1-T_2 \quad \mbox{and} \quad \ddc T_i=\Omega_i^+ -\Omega_i^- \mbox{ for }  i=1,2.
\label{eq_dsh}
\end{eqnarray}
Let $\DSH^p(X)$ denote the space of dsh $(p,p)$-currents on $X$. 
We say that {\it the dsh $(p,p)$-currents
$T_{(n)}$ converge in $\DSH^p(X)$ to $T$} if $T_{(n)}\rightarrow T$ in the sense of currents
and if we can write as in (\ref{eq_dsh})
$$T_{(n)} =T_{1,n}-T_{2,n} \quad \mbox{and} \quad \ddc T_{i,n}=\Omega_{i,n}^+ -\Omega_{i,n}^- 
\mbox{ for }  i=1,2$$
so that the masses of $T_{i,n}$ and $\Omega_{i,n}^\pm$ are locally
uniformly bounded. 
\end{definition}

Dsh currents have been introduced in \cite{DinhSibony1, DinhSibony4, DinhSibony5}.
They are stable under push-forward by holomorphic proper maps and
under pull-back by holomorphic submersions.
A good example of dsh currents to have in mind is the product
$(\varphi_1-\varphi_2) T$ where $\varphi_1$, $\varphi_2$ are  bounded
quasi-psh functions and where $T$ is
a positive closed current. Recall that $\varphi$ is {\it quasi-p.s.h.}
if locally it is a difference of a p.s.h. function with a smooth one.

\begin{lemma} \label{le_pullback_convergence}
Let $(T_n)$ be a sequence of positive smooth $(p,p)$-forms on $\BB_{k'}$ which
converge in $\DSH^p(\BB_{k'})$ to a current $T$. Then $\pi^*(T_n)\wedge [V]$ converge
in $\DSH^{k-l+p}(\BB_k\times \BB_{k'})$ to a current $S$. If $T$ has
no mass on a Borel set $A$ then $S$ has no mass on $\pi^{-1}(A)$.  
\end{lemma}
\begin{proof}
Since the problem is local it is sufficient to consider the case where 0 belongs to
$V$ and where $L$ is a small ball of center 0 as in Lemma \ref{le_pullback_mass}. Let 
$\Pi:\C^k\times \C^{k'}\rightarrow \C^l$ a generic linear
projection. 
Observe that test $(k'+l-p,k'+l-p)$-forms with support in $L$
are generated by forms of type $\varphi\wedge \Pi^*(\psi)$ where 
$\varphi$ is a form of small support in $\BB_k\times\BB_{k'}$ and
$\psi$ is a form of maximal degree on $\C^l$. It follows that slicing
by fibers of $\Pi$
reduces the problem to the case where $l=0$ and $\pi_{|V}$ defines a finite ramified covering
of degree $d_t$ over $\BB_{k'}$. We use here the Lebesgue convergence theorem in order to show that if 
the slices converge pointwise for a fixed $\Pi$ and the mass is
dominated then we have 
the convergence. The problem is to show that there is a unique cluster
point for the sequence $\pi^*(T_n)\wedge [V]$. 

Let $\Sigma\subset V$ be the smallest analytic subset such that $\pi^*(T_n)\wedge[V]$ 
converge to a positive current 
$\tau_0$ on $\BB_k\times\BB_{k'}\setminus\Sigma$. Lemma \ref{le_pullback_mass} applied to $T_n$ and 
$\ddc T_n$ implies that $\tau_0$ and $\ddc \tau_0$ have finite masses
near $\Sigma$. 
Consider the extensions by zero of $\tau_0$ throught 
$\Sigma$ that we 
denote also by $\tau_0$. 
By Theorem \ref{th_extension}, $\ddc\tau_0$ has order zero.
We want to prove that $\Sigma=\varnothing$. Assume that $\Sigma\not = \varnothing$.
Replacing $\BB_k\times \BB_{k'}$ by a small neighbourhood of
some generic (in the Zariski sense) point $a\in\Sigma$ allows to assume that
$\Sigma$ is smooth and $\pi_{|\Sigma}$ is injective.

Now, consider a limit value $\tau$ of $\pi^*(T_n)\wedge [V]$ in 
$\BB_k\times\BB_{k'}$. Then $\tau-\tau_0$ is a positive current with support 
in $\Sigma$. Lemma \ref{le_pullback_mass} applied to $\ddc T_n$ implies that 
$\ddc (\tau-\tau_0)$ has order zero. The support theorem in 
\cite{AlessandriniBassanelli1} extends Federer's theorem on flat
currents \cite{Federer} to currents of order zero with $\ddc$ of order zero. It 
implies that  $\tau-\tau_0$ is a current of $\Sigma$. 
Let $\varphi$ be a test $(k'-p,k'-p)$-form with compact support in $\BB_{k'}$.
We have
\begin{eqnarray}
\lim_{n\rightarrow\infty} \langle \pi^*(T_n)\wedge [V],\pi^*(\varphi)\rangle & = & 
\lim_{n\rightarrow\infty} \langle T_n, \pi_* \big(\pi^*(\varphi)\wedge [V]\big)\rangle \nonumber\\
& = & \lim_{n\rightarrow\infty} \langle T_n, d_t \varphi\rangle  =
d_t\langle T,\varphi\rangle. \label{eq_pullback}
\end{eqnarray}
Hence, $\langle \tau-\tau_0,\pi^*(\varphi)\rangle$ is independent of the choice of $\tau$.
Since $\tau-\tau_0$ is a current on $\Sigma$ and since $\pi_{|\Sigma}$ is injective, the last
identity implies that $\tau-\tau_0$ is independent of the choice of $\tau$.
In other words, $\pi^*(T_n)\wedge [V]$ converge. This contradicts the property that
$\Sigma\not=\varnothing$.

Now we prove that $S$ has no mass on $\pi^{-1}(A)$.
By slicing, one can assume that $\pi_{|V}$ is finite. Let $\Sigma\subset
V$ be an analytic set such that $S$ has no mass on
$\pi^{-1}(A)\setminus \Sigma$. One can choose $\Sigma$ minimal in the sense
that no proper analytic set in $\Sigma$ satisfies the same property.
As above, if $\Sigma$ is not empty, we can assume that it is a smooth
submanifold of $\BB_k\times \BB_{k'}$ and $\pi_{|\Sigma}$ is injective.
Let $\tau_0$ be the restriction of $S$ to
$\BB_k\times\BB_{k'}\setminus\Sigma$ and let $\widetilde \tau_0$ its
extension by zero. Then Remark \ref{rk_extension} implies that
$\widetilde\tau_0$ is dsh and hence 
$S_{|\Sigma}:=S-\widetilde\tau_0$ is a positive dsh current with support in $\Sigma$. By the 
 support theorem \cite{AlessandriniBassanelli1}, this is a current on $\Sigma$.
One deduces from identity (\ref{eq_pullback}) that 
$\langle S_{|\Sigma},\pi^*(\varphi_{|A})\rangle = d_t \langle T,\varphi_{|A\cap\pi(\Sigma)}\rangle$.
Hence, since $T$ has no mass on $A$, the previous integrals vanish. 
The property that $\pi_{|\Sigma}$ is injective implies that 
$S_{|\Sigma}$ has no
mass on $\pi^{-1}(A)$ which contradicts the definition of $\Sigma$.
\end{proof}

\noindent
{\it End of the proof of Theorem \ref{th_main}---}
Let $\pi_1:X\times X'\rightarrow X$ and $\pi_2:X\times X'\rightarrow
X'$ denote the canonical projections. Let $\Gamma$ denote the graph of
$f$ in $X\times X'$.
For $T$ smooth, we have $f^*(T)=(\pi_1)_*\big(\pi_2^*(T)\wedge[\Gamma]\big)$.
We have to show for the general case that $\pi_2^*(T)\wedge [\Gamma]$ is well defined and then define
$f^*(T):=(\pi_1)_*\big( \pi_2^*(T)\wedge [\Gamma]\big)$. Note
that since $\pi_1$ is proper on $\Gamma$, the operator
$(\pi_1)_*$ is well defined on currents supported on $\Gamma$.

On a small open subset $U$ of $X'$, we approximate $T$ by smooth positive closed
(resp. $\ddc$-closed) forms $T_n$ and define 
$$\pi_2^*(T)\wedge [\Gamma]:= \lim_{n\rightarrow\infty} \pi_2^*(T_n)\wedge [\Gamma] \quad \mbox{in } \pi_2^{-1}(U).$$
Lemma \ref{le_pullback_convergence} implies that the limit exists and does not depend on the choice of $T_n$.
This implies also that if $U$ and $U'$ are small open subsets of $X'$, our construction gives two currents 
which coincide in $\pi_2^{-1}(U\cap U')$. 
Hence,  $\pi_2^*(T)\wedge [\Gamma]$
and $f^*(T)$ are  globally well defined.
Since $\pi_2^*(T_n)\wedge [\Gamma]$ are positive and closed (resp. $\ddc$-closed), $\pi_2^*(T)\wedge [\Gamma]$
and $f^*(T)$ are positive and closed (resp. $\ddc$-closed). The
continuity of 
$T\mapsto \pi_2^*(T)\wedge [\Gamma]$, and then the continuity of
$T\mapsto f^*(T)$, follow from Lemma
\ref{le_pullback_convergence}. 
If $T$ has no mass on a set $A\subset X'$, 
Lemma \ref{le_pullback_convergence} shows also that 
 $\pi_2^*(T)\wedge
[\Gamma]$ has no mass on $\pi_2^{-1}(A)\cap\Gamma$. Hence $f^*(T)$ has no
mass on $f^{-1}(A)$. 
\endproof

Using Lemma \ref{le_pullback_convergence}, 
we prove in the same way the following result.

\begin{theorem} \label{th_pullback_dsh}
Let $f$ be as in Theorem \ref{th_main}. Then the pull-back operator 
$$f^*: \DSH^p(X')\rightarrow \DSH^{p}(X)$$ 
is well defined, continuous and commutes with
$\ddc$. Moreover,
if a dsh current $T$ on $X'$ has no mass on a Borel set $A\subset X'$
then $f^*(T)$ has no mass on $f^{-1}(A)$.
\end{theorem}

\begin{remark}\rm
If $g:X'\rightarrow X''$ is another holomorphic map with fibers of
pure dimension $\dim X'-\dim X''$ then the continuity of the pull-back
operator implies that 
$$(g\circ f)^*(T)=f^*(g^*(T))$$
for the classes of current under consideration. Indeed, this identity
holds for smooth forms.
\end{remark}

%%%%%%%%%%%%%%%%%%%%%%%%%%%%%%%%%%%%%%%%%%%%%%

\section{Meromorphic transforms} \label{section_mt}

Meromorphic transforms (MT for short) were considered in \cite{DinhSibony6} in order  
to treat with the same method 
different problems in dynamics and in the study  of 
distribution of varieties (see also \cite{Dinh1, Dinh2}). In this Section we recall the definition of MT
and introduce the pull-back operator on currents associated to a MT.

\begin{definition}\footnote{this definition differs slightly 
    from the definitions given in the previous references}
\rm
A {\it meromorphic transform} 
$F$ of codimension $l$, $0\leq l\leq k-1$, from $X$ onto $X'$ 
is a finite holomorphic chain 
$\Gamma=\sum\Gamma_i$
such that
\begin{enumerate}
\item[-] $\Gamma_i$ is an irreductible analytic 
subset of dimension $k'+l$ of $X\times X'$.
\item[-] $\pi_1$ restricted to each $\Gamma_i$ is proper.   
\item[-] $\pi_2$ restricted to each $\Gamma_i$ is dominant.
\end{enumerate}
\end{definition}

The second item is always
verified when $X'$ is compact.
We do not assume that the $\Gamma_i$ are smooth or distinct. Of course we can write $\Gamma=\sum n_j\Gamma_j'$
where $n_j$ are positive integers and $\Gamma_j'$ are distinct irreducible analytic sets.
Then a generic point in the support $\cup\Gamma_j'$ of $\Gamma$ belongs to a unique $\Gamma_j'$ and $n_j$ is 
called 
{\it the multiplicity} of $\Gamma$ at $x$. 
In what follows we always write $\Gamma$ as $\sum\Gamma_i$. The indices $i$
allow to count the multiplicities.

Define formally $F:=\pi_2\circ(\pi_{1|\Gamma})^{-1}$ and for $A\subset X$ 
and $B\subset X'$
$$F(A):=\pi_2(\pi_1^{-1}(A)\cap\Gamma) \ \mbox{ and }\ F^{-1}(B)=
\pi_1(\pi_2^{-1}(B)\cap\Gamma).$$
So a generic fiber $F^{-1}(x')$ is either empty or an analytic
subset of pure dimension $l$ of $X$.
The sets 
$$I_1:=\{x\in X,\ \dim \pi_1^{-1}(x)\cap\Gamma >k'+l-k\}$$
and
$$I_2:=\{x'\in X',\ \dim \pi_2^{-1}(x')\cap\Gamma >l\}$$
are {\it the first and second indeterminacy sets} of $F$, they are of codimension $\geq 2$.

When $I_2=\varnothing$ we say that $F$ is {\it pure}.
When $\pi_1$ restricted to each $\Gamma_i$ is surjective, we say that $F$ is 
{\it complete}.
A complete MT $F$  of codimension 0 between manifolds of same dimension
is called a {\it meromorphic correspondence}. If moreover $\pi_1$ restricted
to $\Gamma$ is a finite map, then $F$ is called {\it holomorphic 
correspondence}. If generic fibers of $\pi_{1|\Gamma}$ 
contain only one point we 
obtain a dominant meromorphic map from $X$ onto $X'$.

The following proposition extends known results, see e.g. \cite{Meo}.

\begin{proposition} \label{prop_pullback_qpsh}
Let $\varphi$ be a quasi-psh function on $X'$ and
$\Phi$ be a $(p,q)$-form, $k'+l-k\leq p,q \leq k'$, 
whose coefficients are dominated by $|\varphi|$. Let $F$ be a MT as
above. Then the 
$(k-k'-l+p,k-k'-l+q)$-current
$$F^*(\Phi):=(\pi_1)_*\big(\pi_2^*(\Phi)\wedge [\Gamma]\big)$$
is well defined and has locally finite mass.
Let $\Phi_n$ be a sequence of $(p,q)$-forms  whose coefficients are dominated by $\varphi$.
If $\Phi_n\rightarrow\Phi$ in the sense of currents,
then $F^*(\Phi_n)\rightarrow F^*(\Phi)$. If $F$ is complete 
then $F^*(\Phi)$ has $\Lc^1_\loc$ coefficients. 
\end{proposition}
\proof 
Observe that $(\pi_1)_*$ acts continuously on currents
with support in $\Gamma$ since $\pi_1$ restricted to $\Gamma$ is
proper. So, we only have to define $\pi_2^*(\Phi)\wedge [\Gamma]$.

The function $\varphi\circ\pi_2$ is quasi-psh on $X\times X'$. 
Its restriction to any component of $\Gamma$ is not identically $-\infty$. 
Hence, $\varphi\circ\pi_2$ is a quasi-psh function
on $\Gamma$. In particular,  this function is $[\Gamma]$-integrable. 
For quasi-psh functions on singular 
varieties, see \cite{FornaessNarasimhan}. We can also use a desingularization 
$\tau:\widehat\Gamma\rightarrow \Gamma$ and replace $\pi_i$ by $\pi_i\circ\tau$ 
in order to reduce the integrability of $\varphi\circ\pi_2$ to the smooth case. 

We deduce that  $\pi_2^*(\Phi)$ restricted to $\Gamma$ is a form with
coefficients bounded by a quasi-psh function. It follows that
$\pi_2^*(\Phi)\wedge [\Gamma]$ is well defined and has
locally finite mass. 
Hence, $F^*(T)$ is well defined and has locally finite mass. 
Moreover, we have $\pi_2^*(\Phi_n)\wedge [\Gamma]\rightarrow
\pi_2^*(\Phi)\wedge [\Gamma]$ which implies the convergence in the
proposition.

When $F$ is complete,  
$F^*(\Phi)$ has no mass on sets of Lebesgue measure zero. Hence, it has coefficients 
in $\Lc^1_\loc$. Note that when $F$ is not complete we can have
``vertical'' components of $\Gamma$.
\endproof

\begin{corollary} \label{cor_pullback_cohomology}
Let $F$, $p$ and $q$ be as above. Assume that $X$ and $X'$ are compact K{\"a}hler manifolds. Then the operator 
$$F^*:H^{p,q}(X',\C)\rightarrow H^{k-k'-l+p,k-k'-l+q}(X,\C)$$ 
is well defined and is linear.
\end{corollary}
\proof
When $\Phi$ is a smooth closed $(p,q)$-form then $F^*(\Phi)$ is well 
defined and is a closed current of bidegree
$(k-k'-l+p,k-k'-l+q)$ since $F^*$ commutes with the operators
$\partial$ and $\overline\partial$.
It follows that $F^*$ is well defined on the cohomology groups.
\endproof

Assume that $k-k'-l+1\geq 0$. 
This condition is necessary 
so that the pull-back operator on $(1,1)$-currents is meaningful.
The definition of this operator on positive closed $(1,1)$-currents 
under a holomorphic map is classical \cite{Meo} and it can be extended easily to the
case of a MT.

\begin{proposition} \label{prop_pullback1}
Let $T$ be a positive closed current 
of bidegree $(1,1)$ on $X'$.
Then 
$F^*(T)$ is well defined and is positive closed current of bidegree
$(k-k'-l+1,k-k'-l+1)$ which depends continuously on $T$. 
If $X$ and $X'$ are compact K{\"a}hler manifolds, we have $\{F^*(T)\}=
F^*\{T\}$.
\end{proposition}
\proof 
We can write locally $T=\ddc u$ with $u$ psh and 
define 
$$\pi_2^*(T)\wedge [\Gamma]:=\ddc \big((u\circ\pi_2)[\Gamma]\big)
\quad \mbox{and}\quad 
F^*(T):=(\pi_1)_*\big( \pi_2^*(T)\wedge [\Gamma]\big).$$
Since $u\circ\pi_2$ is psh, the current $(u\circ\pi_2)[\Gamma]$
is well defined. Hence, $F^*(T)$ is well defined. 
Since $\ddc$ and $(\pi_1)_*$ commutes, 
the definition is independent of the choice of  $u$.
The current $\pi_2^*(T)\wedge [\Gamma]$ is positive and closed. Hence, so is $F^*(T)$. 
If $T_n\rightarrow T$, we can write locally $T_n=\ddc u_n$ with 
and $u_n\rightarrow u$ in $\Lc^1$. So the continuity is clear.
The assertion on $\{F^*(T)\}$ follows from this continuity  and
a regularization of $T$, see \cite{Demailly2, DinhSibony3}.
\endproof

Using Lemma \ref{le_pullback_convergence}, we easily extends Theorems \ref{th_main} and
\ref{th_pullback_dsh}, and Corollary \ref{cor_main}
to the case of a pure MT. In particular, the following results hold for
dominant meromorphic maps with finite fibers.

\begin{theorem} \label{th_pullback_pure}
Let $F:X\rightarrow X'$ be a pure MT of codimension $l$ between complex manifolds
$X$, $X'$ of dimensions $k$ and $k'$ respectively.
Let $T$ be a positive closed (resp. $\ddc$-closed) $(p,p)$-current on $X'$.  
Then $F^*(T)$ is a positive closed (resp. $\ddc$-closed) current of bidegree 
$(k-k'+p-l,k-k'+p-l)$ which  depends continuously on $T$. Moreover, 
if $T$ has no mass on $F(A)$, with $A\subset X'$ a Borel set, then $F^*(T)$ has no mass on $A$. 
If $X$ and $X'$ are compact K{\"a}hler manifolds, we have $\{F^*(T)\}=F^*\{T\}$.
\end{theorem}

\begin{theorem} \label{th_pullback_pure_dsh}
Let $F$ and $A$ be as in Theorem \ref{th_pullback_pure}. Then the pull-back operator 
$$F^*: \DSH^p(X')\rightarrow \DSH^{k-k'+p-l}(X)$$ 
is well defined, continuous and commutes with
$\ddc$. 
Moreover, if a dsh current $T$ on $X'$ has no mass on $F(A)$,
then $F^*(T)$ has no mass on $A$.
\end{theorem}

%%%%%%%%%%%%%%%%%%%%%%%%%%%%%%%%%%%%%%%%%%%%%%%%%%%%%%%%%

\section{The case of compact K\"ahler manifolds} \label{section_karler}

In general it is not possible to define the pull-back of a current under a holomorphic map, in 
a consistent way. A simple example in \cite{Meo} shows that the mass of the
pull-back could be infinite around exceptional fibers. However the situation in the compact case is different.

From now on $(X,\omega)$ and $(X',\omega')$ denote compact K\"ahler
manifolds of dimension $k$ and $k'$ respectively.
Consider a MT $F:X\rightarrow X'$ of codimension $l$ of graph $\Gamma$. Let $\Cc$ be the
analytic subset of $\Gamma$
defined by the property that $\pi_2$ restricted to 
$\Gamma\setminus\Cc$ has locally only empty fibers or fibers of dimension pure
$l$.  
Let $\pi$ denote the restriction of $\pi_2$ to $\Gamma\setminus\Cc$.
If $T$ is a positive closed (resp. $\ddc$-closed) 
current of bidegree $(p,p)$, $k'+l-k\leq p\leq k'$, on $X'$, then by
Theorem \ref{th_pullback_pure}, the current 
$\pi^*(T)$ is well defined on $\Gamma\setminus\Cc$. 
We will show 
that $\pi^*(T)$ has finite mass; this allows us to extend $\pi^*(T)$ through $\Cc$. 
Denote by $(\pi_{2|\Gamma})^\star (T)$ the trivial extension of $\pi^* (T)$. 
By Theorem \ref{th_extension}, 
$\ddc (\pi_{2|\Gamma})^\star (T) \leq 0$. In our situation, since $X'$ is compact
Stokes formula implies that for any closed form $\Omega$ of right bidegree 
$$\langle \ddc (\pi_{2|\Gamma})^\star (T), \Omega \rangle = \langle 
(\pi_{2|\Gamma})^\star (T), \ddc \Omega \rangle = 0.$$
In particular, the mass of $\ddc (\pi_{2|\Gamma})^\star (T)$ is zero, then
$(\pi_{2|\Gamma})^\star (T)$ is $\ddc$-closed.
We call  {\it the strict transform} of $T$ by $F$ the current
$$F^\star (T):=(\pi_1)_*(\pi_{2|\Gamma})^\star (T).$$

\begin{proposition}[see also \cite{DinhSibony2, DinhSibony3}] \label{prop_strict_transform}
Let $F$ and  $T$ be as above. Then 
$F^\star (T)$ is positive closed (resp. $\ddc$-closed). 
Moreover, there exists a constant $c>0$ independent of $T$
such that $\|F^\star (T)\|\leq c\|T\|$.
The operator $T\mapsto F^\star (T)$ is lower semi-continuous in the sense that 
if $T_n\rightarrow T$ then every cluster value $\tau$ of $F^\star (T_n)$ 
satisfies $\tau\geq F^\star(T)$. If $F$ is complete and if $T$ has no mass on
analytic sets then $F^\star(T)$ has no mass on analytic sets. 
\end{proposition}
\begin{proof} We first prove that $(\pi_{2|\Gamma})^\star (T)$ has finite mass.
By \cite{DinhSibony3}, there are positive closed (resp. $\ddc$-closed)
smooth forms $T_n^\pm$ with cohomology classes bounded by $c'\|T\|$,
$c'>0$, such that $T_n^\pm\rightarrow T^\pm$ and $T^+-T^-=T$.
We have 
$$\|(\pi_{2|\Gamma})^\star (T)\|\leq \|(\pi_{2|\Gamma})^\star (T^+)\|=
\lim_{n\rightarrow\infty} \|(\pi_{2|\Gamma})^\star (T_n^+)\|\leq \lim_{n\rightarrow\infty} 
\|(\pi_2)^*(T_n^+)\wedge [\Gamma]\|.$$ 
Since 
the cohomological classes $\{T_n^+\}$
of $T_n^+$ are uniformly bounded with respect to $n$, $\|(\pi_2)^*(T_n^+)\wedge [\Gamma]\|$, which can be computed 
cohomologically, are uniformly 
bounded. More precisely, we have $\|(\pi_2)^*(T_n^+)\wedge [\Gamma]\|\leq c''\|T\|$, $c''>0$, for
$n$ large enough. Hence $F^\star (T)$ is well defined and 
$\|F^\star (T)\|\leq c\|T\|$ with $c>0$ independent of $T$. 

By definition, $(\pi_{2|\Gamma})^\star (T) =\lim (\pi_{2|\Gamma})^\star (T_n)$ 
on $\Gamma\setminus \Cc$.
Since $(\pi_{2|\Gamma})^\star (T)$ has no mass on $\Cc$, it follows that 
$(\pi_{2|\Gamma})^\star (T)$
is smaller than any cluster value of $(\pi_{2|\Gamma})^\star (T_n)$. 
Hence $\tau\geq F^\star (T)$. The last statement follows from Theorem 
\ref{th_pullback_pure}.
\end{proof}

\begin{definition} \rm
If $\{F^\star (T)\} = F^* \{T\}$ we say
that $F^*(T)$ is {\it well defined} and we write $F^*(T):=F^\star
(T)$.
We call $F^*(T)$ {\it the total transform} of $T$. 
\end{definition}
For currents in a projective space, this definition has been 
used in order to study the 
dynamics of birational maps \cite{DinhSibony5} 
(see also \cite{AlessandriniBassanelli2}). The following result
justifies the previous definition.

\begin{proposition} The operator $F^*$ is continuous in the following sense. 
Let $T_n$ and $T$  be
positive closed ($\ddc$-closed) currents such that $T_n\rightarrow T$. Assume that $F^*(T_n)$ and 
$F^*(T)$ are well defined in the above sense. Then 
$F^*(T_n)\rightarrow F^* (T)$. 
\end{proposition}
\begin{proof} As in Proposition \ref{prop_strict_transform}, $\|F^*(T_n)\|$ is bounded
  uniformly on $n$. We can assume that 
$F^*(T_n)$ converge 
to a current $\tau$.
Proposition \ref{prop_strict_transform} 
implies that $\tau\geq F^*(T)$. On the other hand, since $F^*$ acts continuously on cohomology
groups, we have 
$$\{\tau\}=\lim\{F^*(T_n)\}=\lim F^*\{T_n\}=F^*\{T\}=\{F^*(T)\}.$$
It follows that $\|\tau\|=\|F^*(T)\|$ and hence $\tau=F^*(T)$.
\end{proof}

Now we study the case of $\ddc$-closed $(1,1)$-currents.

\begin{theorem} \label{th_pullback_11}
Let $F:X\rightarrow X'$ be a MT between compact K{\"a}hler manifolds. If $T$ is a positive 
$\ddc$-closed $(1,1)$-current on $X'$. Then there is a unique
$\ddc$-closed extension $F^*(T)$ of $F^\star(T)$ such that 
$T\mapsto F^*(T)$ is continuous for the weak topology on currents. 
Moreover, we have $\{F^*(T)\}=F^*\{T\}$.
\end{theorem}

For the proof we need the following fact.

\begin{proposition} \label{prop_class_exceptional}
Let $f:X\rightarrow X'$ be a holomorphic surjective map between
compact K\"ahler manifolds of dimension $k$
and $k'$. 
Let $I_2$ be the set of $x'\in X'$ such that $\dim f^{-1}(x')>k-k'$. 
Then the components of codimension $1$ of $f^{-1}(I_2)$ are cohomologically independent.
\end{proposition}
\proof
By Stein's factorization theorem \cite[E.G.A III 4.3.3]{Grothendieck}, there exist a normal space $Y$,
a holomorphic map $h:X\rightarrow Y$ and a finite morphism $g:Y\rightarrow X'$ such that $f=g\circ h$. Moreover
generic fibers of $h$ are connected. We can replace $X'$ by $Y$ and
assume that generic fibers of $f$ are connected, but $X'$ may have
singularities \footnote{we thank F. Campana who told us this argument}.

Let $E_i$ be components of codimension 1 of $f^{-1}(I_2)$. We have to
show that the classes $\{E_i\}$ are linearly independent. Let 
$c_i\in\R$
such that the class of $\sum c_i E_i$ vanishes. Then there is an integrable
function $u$ 
such that $\ddc u=\sum c_i [E_i]$.
This function is pluriharmonic out of $f^{-1}(I_2)$.

Let $\Omega$ be a closed smooth $(k-k',k-k')$-form on $X$. Then 
$f_*(u\Omega)$ is $\ddc$-closed. It follows that the current
$f_*(u\Omega)$ is equal to a function which is pluriharmonic out of $I_2$. 
Since $I_2$ has codimension $\geq 2$ this
function can be extended to a pluriharmonic function on the compact
space $X'$ and hence should be constant.

We show that  $u$ is 
constant on each generic fiber of $f$. The case $k=k'$ is clear. One
consider the case $k>k'$. In a neighbourhood of a generic
point $x\in X$ one can find a coordinate system $(z,z')$ so that $x=0$
and 
$f(z,z')=z'$. If $\varphi$ is a smooth function with support in this
neighbourhood, then $\Omega:=\ddc\varphi \wedge (\ddc\|z\|^2)^{k-k'-1}$ is a
closed form. Since $f_*(u\Omega)$ is a constant function vanishing
near $I_2$, we have $f_*(u\Omega)=0$. It follows that 
$$\int u(z,0)
\ddc\varphi(z,0)  \wedge (\ddc\|z\|^2)^{k-k'-1}=0$$ 
for every $\varphi$. Hence $u(z,0)$ is constant. Since generic fibers of $f$
are connected, $u$ is constant on each generic fiber of $f$. Then
there is a function $u'$ on $X'$ such that $u=u'\circ f$ almost
everywhere. 

Now consider a strictly positive closed form $\Omega$. 
We have $f_*(u\Omega)=u'f_*(\Omega)$. On the other hand, since $f_*(\Omega)$
is a closed $(0,0)$-current, it is given by a constant function. 
This function is not zero because $\Omega$ is strictly positive.
Finally, the fact that $f_*(u\Omega)$ is given by a constant function implies
that $u'$ is constant and then $u$ is constant. Consequently, $c_i=0$ for
every $i$.
\endproof
\noindent
{\it Proof of Theorem \ref{th_pullback_11}.} 
We can assume that the graph $\Gamma$ of $F$ is irreducible
and that $\Gamma$ is smooth. Otherwise, we consider a blow-up
$\tau:\widetilde\Gamma\rightarrow\Gamma$ and use $\pi_i\circ \tau$
instead of $\pi_{i|\Gamma}$.

By \cite{DinhSibony3}, $T$ is a difference of currents which can be
approximated by smooth positive
$\ddc$-closed $(1,1)$-forms. 
We can assume that $T$ is the limit of smooth positive $\ddc$-closed 
$(1,1)$-forms $T_n$. 
Recall that $F^*$ is well defined on smooth forms. Then the uniqueness
in the theorem is clear.

The masses of $(\pi_{2|\Gamma})^*(T_n)$ are computed cohomologically;
they are bounded uniformly on $n$.
We only have to show that $(\pi_{2|\Gamma})^*(T_n)$ converge. If $\tau$ is a limit value then 
$\{\tau\}=(\pi_{2|\Gamma})^*\{T\}$ and $\tau-(\pi_{2|\Gamma})^\star(T)$ is a pluriharmonic $(1,1)$-current
with support in $\pi_2^{-1}(I_2)\cap\Gamma$. If the $E_i$ are components of codimension 1 of
$\pi_2^{-1}(I_2)\cap\Gamma$ 
then there are real numbers
$c_i$ such that
$\tau=  (\pi_{2|\Gamma})^\star(T)+\sum c_i[E_i]$. Proposition \ref{prop_class_exceptional} implies that $c_i$ are
uniquely determined by the identity $\{\tau\}=\{(\pi_{2|\Gamma})^*T\}+\sum c_i \{[E_i]\}$. The theorem follows.
\hfill $\square$

%%%%%%%%%%%%%%%%%%%%%%%%%%%%%%%%%%%%%%%%%%%%%%%%%%%%

Tien-Cuong Dinh \hfill Nessim Sibony\\
Institut de Math{\'e}matique de Jussieu \hfill  Math{\'e}matique - B{\^a}timent 425 \\
Plateau 7D, Analyse Complexe \hfill  UMR 8628\\
175 rue du Chevaleret \hfill  Universit{\'e} Paris-Sud \\
75013 Paris, France \hfill 91405 Orsay, France \\
{\tt dinh@math.jussieu.fr} \hfill {\tt nessim.sibony@math.u-psud.fr} \\


\small
\begin{thebibliography}{11}


\addcontentsline{toc}{section}{References}


\bibitem{AlessandriniBassanelli1}
Alessandrini, L. and Bassanelli, G.
Plurisubharmonic currents and their extension across analytic subsets,
{\it Forum Math.}, {\bf 5} (1993), no. 6, 577-602.

\bibitem{AlessandriniBassanelli2}
Alessandrini, L. and Bassanelli, G. 
Transforms of currents by modifications and 1-convex manifolds,
{\it  Osaka J. Math.}, {\bf 40} (2003), no. 3, 717-740.


\bibitem{Bassanelli}
Bassanelli, G. A cut-off theorem for plurisubharmonic currents,  
\textit{Forum Math.},  \textbf{6}  (1994),  no. 5, 567--595.

\bibitem{ChernLevineNirenberg}
Chern, S.S. Levine, H.I. and Nirenberg, L. 
Intrinsic norms on a complex manifold, in 
Global Analysis (Papers in Honor of K. Kodaira), Univ. Tokyo Press, 1969,  119-139. 

\bibitem{DabbekElkhadharaElMir}
Dabbek, K., Elkhadhra, F. and El Mir, H. Extension of plurisubharmonic currents, {\it Math. Z.}, 
{\bf 245} (2003), 455-481.


\bibitem{Demailly2}
Demailly, J.-P. {\it Complex analytic geometry}, available at \\
www.fourier.ujf-grenoble.fr/$\sim$demailly. 

\bibitem{Dinh1}
Dinh, T.-C. Distribution des pr{\'e}images et des points
p{\'e}riodiques d'une correspondance polynomiale, 
{\it Bull. Soc. Math. France}, {\bf 133} (2005), no. 3, 363-394. 

\bibitem{Dinh2}
Dinh, T.-C. Suites d'applications m{\'e}romorphes multivalu{\'e}es et 
courants laminaires, \textit{J. Geom. Analysis}, {\bf 15} (2005), no. 2, 207-227.

\bibitem{DinhSibony1}
Dinh, T.-C. and Sibony, N. 
Dynamique des applications d'allure
polynomiale, {\it J. Math. Pures Appl.}, {\bf 82} (2003),
367-423. 

\bibitem{DinhSibony2}
Dinh, T.-C. and Sibony, N. Une borne sup{\'e}rieure pour
l'entropie topologique  d'une application rationnelle, {\it Ann. of
  Math.(2)} {\bf 161} (2005), no. 3, 1637-1644.

\bibitem{DinhSibony3}
Dinh, T.-C. and Sibony, N. Regularization of currents and entropy, 
{\it Ann. Sci. Ecole Norm. Sup.},  {\bf 37} (2004), 959-971. 


\bibitem{DinhSibony4}
Dinh, T.-C. and Sibony, N. Green currents for automorphisms of 
compact K{\"a}hler manifolds, \textit{J. Amer. Math. Soc.},  {\bf 18}
(2004), 291-312. 

\bibitem{DinhSibony5}
Dinh, T.-C. and Sibony, N. Dynamics of regular birational maps in $\mathbb{P}^k$, 
{\it J. Funct. Anal.},  {\bf 222} (2005), 202-216.


\bibitem{DinhSibony6}
Dinh, T.-C. and  Sibony, N. Distribution de valeurs d'une
suite de transformations m{\'e}romorphes et applications, {\it
  Comment. Math. Helv.}, {\bf 81} (2006), 221-258.


\bibitem{ElMir}
El Mir, H.
Sur le prolongement des courants positifs ferm{\'e}s,
{\it Acta Math.}, {\bf 153} (1984), no. 1-2, 1-45.

\bibitem{Federer}
Federer, H. {\it Geometric Measure Theory}, 
New York, Springer-Verlag, 1969.

\bibitem{FornaessNarasimhan}
Forn\ae ss,  J.-E. and Narasimhan, R.
The Levi problem on complex spaces with singularities,
{\it Math. Ann.}, {\bf 248} (1980), no. 1, 47-72.


\bibitem{Grothendieck}
Grothendieck, A. {\it EGA}, S\'eminaire du Bois-Marie IHES, 1967.


\bibitem{Meo}
M{\'e}o, M. Image inverse d'un courant positif ferm{\'e} par une
application surjective, \textit{C.R.A.S.}, \textbf{322} (1996),
1141-1144. 


\bibitem{Sibony2}
Sibony, N. Quelques probl{\`e}mes de 
prolongement de courants en analyse complexe, 
\textit{Duke Math. J.}, \textbf{52} (1985),  no.1, 157--197.

\bibitem{Sibony1}
Sibony, N. Dynamique des applications rationnelles de
$\mathbb{P}^k$, \textit{Panoramas et Synth{\`e}ses}, {\bf 8} (1999), 97-185.


\bibitem{Skoda1}
Skoda, H. Prolongement des courants positifs, ferm{\'e}s de
masse finie, \textit{Invent. Math.}, \textbf{66} (1982), 361-376.





\end{thebibliography}
\normalsize
\end{document}